

\baselineskip=14pt
\parskip=10pt

\magnification=\magstephalf

\def\1{{\overline{1}}}
\def\2{{\overline{2}}}
\parindent=0pt
\overfullrule=0in

\def\frac#1#2{{#1 \over #2}}
\centerline
{\bf  Answers to Some Questions about  Explicit Sinkhorn Limits posed by Mel Nathanson } 
\bigskip
\centerline
{\it Shalosh B. EKHAD and Doron ZEILBERGER}
\bigskip

{\bf Preface}

At the Jan. 2018 Joint Mathematics Meetings, Avi Wigderson 
gave a series of three  fascinating lectures [W], whose starting point was
the {\it Sinkhorn algorithm}. One of the people in the audience was
Mel Nathanson, and this lead him to write two  papers [N1][N2] inspired by this algorithm.

A square matrix is {\it row-stochastic} if all its rows add-up to $1$ and is
{\it column-stochastic} if all its columns add-up to $1$. It is {\it doubly-stochastic} if
it is both row- and column- stochastic.

Consider the following question.

{\it
``Given a square  matrix, $A$,  with positive entries, can you find diagonal matrices $X$ and $Y$, 
and a doubly-stochastic matrix, $S$, such $S=XAY$?}''

{\bf  The Sinkhorn algorithm} gives, very fast, an {\it approximate} answer, as follows.
Let $R(A)$ be the operation that inputs a matrix $A$ with positive entries and outputs the
row-stochastic matrix obtained by normalizing each row, i.e. dividing each row by its sum.
Analogously, let $C(A)$, be the operation that inputs a matrix $A$ with positive entries and outputs the
column-stochastic matrix obtained by normalizing each column, i.e. dividing each column by its sum.

The Sinkhorn algorithm proceeds by {\it alternating} these two `correction' operations. Surprisingly, after few iterations
you get something that is {\it approximately} doubly-stochastic.

If $A$ is also symmetric,  then one can take $X=Y$, and then one is looking for the unique diagonal matrix, $X$, and
for the unique symmetric doubly-stochastic matrix, $S$, such that $S= XAX$.

{\bf Mel Nathanson's questions}

Nathanson wondered if one can find explicit expressions, in terms of the
entries of $A$, for the entries of the Sinkhorn limit, $S$. 
He also commented that one should be able to do it using {\it Gr\"obner bases}.  He also  wondered
whether there exist matrices for which the Sinkhorn algorithm
terminates in a {\it finite} number of steps, and settled the question for the $2 \times 2$ case.

{\bf The Maple package SINKHORN.txt}

One of us (DZ) wrote a Maple package, {\tt SINKHORN.txt}, available from the front of the present article

{\tt http://sites.math.rutgers.edu/\~{}zeilberg/mamarim/mamarimhtml/sinkhorn.html} \quad ,

that lead to the solutions, by the other author (SBE), of some of Nathanson's questions.
Following Nathanson's advice we used Gr\"obner bases (the Buchberger algorithm).

{\bf Answers to some of Nathanson's questions}

The following theorem completely answers  the {\it central problem} (problem 1, p. 26) in Nathanson's article [N2].

{\bf Theorem 1.} Let 
$$
A= \left ( \matrix{a_{11} & a_{12} & a_{13} \cr
                   a_{12} & a_{22} & a_{23} \cr
                   a_{13} & a_{23} & a_{33}}
     \right )
$$
be the {\it generic}, (`symbolic'), $3 \times 3$ {\bf symmetric} matrix, with {\bf positive coefficients}.
Its Sinkhorn limit, let's call it $S$,  is a certain symmetric $3 \times 3$ {\bf doubly-stochastic} matrix whose $(1,1)$ entry, $s_{11}$, is given by
$$
s_{11} = a_{11} z \quad,
$$

where $z$ is the positive root of the {\bf quartic equation}
$$
( {a_{{1,1}}}^{4}{a_{{2,2}}}^{2}{a_{{3,3}}}^{2}-{a_{{1,1}}}^{4}a_{{2,2}}{a_{{2,3}}}^{2}a_{{3,3}}-2\,{a_{{1,1}}}^{3}{a_{{1,2}}}^{2}a_{{2,2}}{a_{{3,3}}}
^{2}
$$
$$
+{a_{{1,1}}}^{3}{a_{{1,2}}}^{2}{a_{{2,3}}}^{2}a_{{3,3}}+2\,{a_{{1,1}}}^{3}a_{{1,2}}a_{{1,3}}a_{{2,2}}a_{{2,3}}a_{{3,3}}-2\,{a_{{1,1}}}^{3}{a_{{1,3}}}^{2
}{a_{{2,2}}}^{2}a_{{3,3}}+{a_{{1,1}}}^{3}{a_{{1,3}}}^{2}a_{{2,2}}{a_{{2,3}}}^{2}
$$
$$
+{a_{{1,1}}}^{2}{a_{{1,2}}}^{4}{a_{{3,3}}}^{2}-2\,{a_{{1,1}}}^{2}{a_{{1,2}}}
^{3}a_{{1,3}}a_{{2,3}}a_{{3,3}}+3\,{a_{{1,1}}}^{2}{a_{{1,2}}}^{2}{a_{{1,3}}}^{2}a_{{2,2}}a_{{3,3}}-{a_{{1,1}}}^{2}{a_{{1,2}}}^{2}{a_{{1,3}}}^{2}{a_{{2,3}}}^
{2}
$$
$$
-2\,{a_{{1,1}}}^{2}a_{{1,2}}{a_{{1,3}}}^{3}a_{{2,2}}a_{{2,3}}+{a_{{1,1}}}^{2}{a_{{1,3}}}^{4}{a_{{2,2}}}^{2}-a_{{1,1}}{a_{{1,2}}}^{4}{a_{{1,3}}}^{2}a_{{3,
3}}+2\,a_{{1,1}}{a_{{1,2}}}^{3}{a_{{1,3}}}^{3}a_{{2,3}}-a_{{1,1}}{a_{{1,2}}}^{2}{a_{{1,3}}}^{4}a_{{2,2}} ) \, {z}^{4}
$$
$$
+ ( -4\,{a_{{1,1}}}^{3}{a_{{2,2}}}^{2}{a_{{3,3}}}^{2}+4\,{a_{{1,1}}}^{3}a_{{2,2}}{a_{{2,3}}}^{2}a_{{3,3}}+4\,{a_{{1,1}}}^{
2}{a_{{1,2}}}^{2}a_{{2,2}}{a_{{3,3}}}^{2}
-3\,{a_{{1,1}}}^{2}{a_{{1,2}}}^{2}{a_{{2,3}}}^{2}a_{{3,3}}-2\,{a_{{1,1}}}^{2}a_{{1,2}}a_{{1,3}}a_{{2,2}}a_{{2,3}}a_
{{3,3}}
$$
$$
+4\,{a_{{1,1}}}^{2}{a_{{1,3}}}^{2}{a_{{2,2}}}^{2}a_{{3,3}}-3\,{a_{{1,1}}}^{2}{a_{{1,3}}}^{2}a_{{2,2}}{a_{{2,3}}}^{2}-2\,a_{{1,1}}{a_{{1,2}}}^{2}{a_{{
1,3}}}^{2}a_{{2,2}}a_{{3,3}}+2\,a_{{1,1}}{a_{{1,2}}}^{2}{a_{{1,3}}}^{2}{a_{{2,3}}}^{2}-{a_{{1,2}}}^{4}{a_{{1,3}}}^{2}a_{{3,3}}
$$
$$
+2\,{a_{{1,2}}}^{3}{a_{{1,3}}}
^{3}a_{{2,3}}-{a_{{1,2}}}^{2}{a_{{1,3}}}^{4}a_{{2,2}} ) \, {z}^{3}
$$
$$
+ ( 6\,{a_{{1,1}}}^{2}{a_{{2,2}}}^{2}{a_{{3,3}}}^{2}-6\,{a_{{1,1}}}^{2}a_{{2,2}}{a_{{2,3}}}^{2}a_{{3,3}}-2\,a_{{1,1}}{a_{{1,2}}}^{2}a_{{2,2}}{a_{{3,3}}}^{2}
+3\,a_{{1,1}}{a_{{1,2}}}^{2}{a_{{2,3}}}^{2}a_{{3,3}}
$$
$$
-2\,a_{{1,1}}a_{{1,2}}a_{{1,3}}a_{{2,2}}a_{{2,3}}a_{{3,3}}-2\,a_{{1,1}}{a_{{1,3}}}^{2}{a_{{2,2}}}^{2}a_{{3,3}}+3\,a_{{1,1}}{a_{{1,3}}}^{2}a_{{2,2}}{a_{{2,3}}}^{2}
+2\,{a_{{1,2}}}^{3}a_{{1,3}}a_{{2,3}}a_{{3,3}}
$$
$$
-3\,{a_{{1,2}}}^{2}{a_{{1,3}}}^{2}a_{{2,2}}a_{{3,3}}-{a_{{1,2}}}^{2}{a_{{1,3}}}^{2}{a_{{2,3}}}^{2}+2\,a_{{1,2}}{a_{{1,3}}}^{3}a_{{2,2}}a_{{2,3}} ) \, {z}^{2}
$$
$$
+ ( -4\,a_{{1,1}}{a_{{2,2}}}^{2}{a_{{3,3}}}^{2}+4\,a_{{1,1}}a_{{2,2}}{a_{{2,3}}}^{2}a_{{3,3}}-{a_{{1,2}}}^{2}{a_{{2,3}}}^{2}a
_{{3,3}}+2\,a_{{1,2}}a_{{1,3}}a_{{2,2}}a_{{2,3}}a_{{3,3}}-{a_{{1,3}}}^{2}a_{{2,2}}{a_{{2,3}}}^{2} ) \, z
$$
$$
+{a_{{2,2}}}^{2}{a_{{3,3}}}^{2}-a_{{2,2}}{a_{{2,3}}}^{2}a_{{3,3}} \, = \, 0 \quad .
$$

Furthermore the diagonal matrix $X$, such that $S=XAX$ has its $(1,1)$ entry, $x_{11}$, given explicitly by
$$
x_{11}= \sqrt{z} \quad .
$$

The other five entries of the $3 \times 3$ symmetric matrix $S$, and the other two non-zero entries of the  $3 \times 3$ diagonal matrix $X$ are
too long to be presented here, but are readily available (free of charge, and no advertisements!) from the following url

{\tt http://sites.math.rutgers.edu/\~{}zeilberg/tokhniot/oSINKHORN3.txt} \quad .

{\bf Comment}: Of course `explicit' is in the eyes of the beholder, and some people may argue that Sinkhorn's
algorithm that produces (extremely fast!) the desired doubly-stochastic matrix $S$ to any desired accuracy
is explicit enough. But to pure mathematicians it only gives `approximations'. Our solution is as explicit as
it can get, even if you insist that the entries are `solvable by radicals', since $z$ satisfies
a certain explicit {\it quartic} equation, with coefficients that are polynomials in the six entries of 
$A$.

Since the general case is so complicated, Nathanson [N2] (problem 2, p. 26) also asked for the Sinkhorn matrices
of two special cases. The next theorem answers the first part of problem 2.

{\bf Theorem 2.} Let  $K$ and $L$ be arbitrary positive numbers, and let
$$
A= \left ( \matrix{K & 1 & 1 \cr
                   1 & L & 1 \cr
                   1 & 1 & 1}
     \right ) \quad .
$$
Its Sinkhorn limit, let's call it $S$,  is a certain symmetric $3 \times 3$ {\bf doubly-stochastic} matrix whose $(1,1)$ entry, $s_{11}$, is given by
$$
s_{11} = K z \quad,
$$
where $z$ is the positive root of the quartic equation
$$
L+ \left( -4\,LK+1 \right) z+ \left( 6\,L{K}^{2}-2\,LK-3\,K-1 \right) {
z}^{2}- \left( K-1 \right)  \left( 4\,L{K}^{2}-3\,K-1 \right) {z}^{3}+K
 \left( K-1 \right) ^{2} \left( LK-1 \right) {z}^{4}=0 \quad .
$$

Furthermore the diagonal matrix $X$, such that $S=XAX$ has its $(1,1)$ entry, $x_{11}$, given explicitly by
$$
x_{11}= \sqrt{z} \quad .
$$

The other five entries of the $3 \times 3$ symmetric matrix $S$, and the other two non-zero entries of the  $3 \times 3$ diagonal matrix $X$ are
available here: 

{\tt http://sites.math.rutgers.edu/\~{} zeilberg/tokhniot/oSINKHORN4.txt } \quad .

The next theorem answers the second  part of problem 2 of [N2].

{\bf Theorem 3.} Let  $K$, $L$ and $M$ be arbitrary positive numbers, and let
$$
A= \left ( \matrix{K & 1 & 1 \cr
                   1 & L & 1 \cr
                   1 & 1 & M}
     \right )
$$
Its Sinkhorn limit, let's call it $S$,  is a certain symmetric $3 \times 3$ {\bf doubly-stochastic} matrix whose $(1,1)$ entry, $s_{11}$, is given by
$$
s_{11} = K z \quad,
$$
where $z$ is the positive root of the quartic equation
$$
LM ( LM-1 ) + ( -4\,K{L}^{2}{M}^{2}+4\,KML+2\,LM-L-M ) z
$$
$$
+ ( 6\,{K}^{2}{L}^{2}{M}^{2}-6\,{K}^{2}LM-2\,KM{L}^{2}-2
\,KL{M}^{2}-2\,KML+3\,LK+3\,KM-3\,LM+2\,L+2\,M-1 ) {z}^{2}+
$$
$$
( -4\,{K}^{3}{L}^{2}{M}^{2}+4\,{K}^{3}LM+4\,{K}^{2}{L}^{2}M+4\,{K}^{2}L{M}^{2}-2\,{K}^{2}LM-3\,{K}^{2}L
$$
$$
-3\,{K}^{2}M-2\,KML+2\,K-L-M+2
 ) {z}^{3}
$$
$$
+ ( KM-1 ) K ( LK-1 )  ( KML
-K-L-M+2 ) {z}^{4}=0 \quad .
$$

Furthermore the diagonal matrix $X$, such that $S=XAX$ has its $(1,1)$ entry, $x_{11}$, given explicitly by
$$
x_{11}= \sqrt{z} \quad .
$$

The other five entries of the $3 \times 3$ symmetric matrix $S$, and the other two entries of the diagonal matrix $X$ are
available here:

{\tt http://sites.math.rutgers.edu/\~{} zeilberg/tokhniot/oSINKHORN5.txt } \quad .

The next fact answers, in the {\bf affirmative}, problem 5 (p. 27) in [N2].

{\bf Fact 4}: The matrix
$$
A \, = \,
\left ( \matrix{ \frac{1}{5} & \frac{1}{5} & \frac{3}{5} \cr
                             &             &             \cr
               \frac{2}{5} & \frac{1}{5} & \frac{2}{5} \cr
                             &             &             \cr
               \frac{3}{5} & \frac{1}{5} & \frac{1}{5} }
\right ) 
$$
is row-stochastic (check!), but not column-stochastic (check!), but applying column-scaling to it yields the matrix
$$
\left ( \matrix{ \frac{1}{6} & \frac{1}{3} & \frac{1}{2} \cr
                             &             &             \cr
               \frac{1}{3} & \frac{1}{3} & \frac{1}{3} \cr
                             &             &             \cr
               \frac{1}{2} & \frac{1}{3} & \frac{1}{6} }
\right ) 
$$
that {\bf is} doubly-stochastic (check!). 

By multiplying the first row of $A$ by $10$, the second row by $5$ and the third row by $15$ we get the matrix
$$
M \, = \,
\left ( \matrix{ 2 & 2 & 6 \cr
                 2 & 1 & 2 \cr
                 9 & 3 & 3}
\right )  \quad ,
$$
that achieves its Sinkhorn limit after only {\bf two} steps (or one double step).
In other words $M$ is not doubly-stochastic but $C(R(M))$ is.

Using procedure {\tt MelNprob5(T,var)} in the Maple package {\tt SINKHORN.txt}, one can concoct many other such examples.

{\bf On the $3 \times 3$ matrix of all $1$s except for the $(1,1)$ entry}

In section 13 of [N2], the following matrix is discussed
$$
A(r):=\left ( \matrix{ \frac{r(r+1)}{2} & 1 & 1 \cr
                             &             &             \cr
                              1  & 1 & 1 \cr
                             &             &             \cr
                              1  & 1 & 1} 
 \right )\quad .
$$
Its Sinkhorn limit , let's call it $S(r)$ is:
$$
S(r)=
\left(
\matrix{
\frac{r}{r+2} & \frac{1}{r+2} & \frac{1}{r+2} \cr
                             &             &             \cr
\frac{1}{r+2} & \frac{r+1}{2(r+2)} & \frac{r+1}{2(r+2)} \cr
                             &             &             \cr
\frac{1}{r+2} & \frac{r+1}{2(r+2)} & \frac{r+1}{2(r+2)}
       }
 \right ) \quad .
$$

The diagonal matrix, $X(r)$, such that $X(r)A(r)X(r)=S(r)$ is
$$
X(r)=
\sqrt{\frac{2}{(r+1)(r+2)}} \,\,
\left (
\matrix{1 & 0 & 0 \cr
                             &             &             \cr
        0 & \frac{r+1}{2} & 0 \cr
                             &             &             \cr
        0 &   0           &\frac{r+1}{2}}
\right ) \quad .
$$

It is asked in [N2] whether the Sinkhorn algorithm applied to $A(r)$ 
can terminate after a finite number of steps.
This is unlikely for the following reason. We Use procedure {\tt MelNsec13(r,k)} in our Maple package,
It inputs a symbol {\tt r}, and a positive integer {\tt k}, and outputs the difference between
the  sums of the first and second rows when row-scaling followed by column-scaling is applied
{\tt k} times. This is a necessary condition for being doubly-stochastic.
By trying out  {\tt MelNsec13(r,k)} for {\tt k} from $1$ to $6$ it appears that the numerator is always
$3\,((r+2)(r-1))^{2k}$, hence only vanishes when $r=1$ or $r=-2$ producing the all $1$-matrix.
The fact that this holds for all $k$ could presumably proved rigorously by mathematical induction.

{\bf What about Larger sizes?}

Theorem 1 was obtained via procedure {\tt ExacGS} in our Maple package {\tt SINKHORN.txt} .
It would be too much for Maple (and probably also for SINGULAR and even for MAGMA) to do the analogous
theorem for a {\it generic}, symbolic symmetric $n \times n$ matrix for $n \geq 4$. But it does a good job,
for {\it numerical} matrices, finding the {\it exact} Sinkhorn limits in terms of algebraic numbers.

If one had a sufficiently large computer, one would be able to state the analog of Theorem 1 for
$n \times n$ matrices, for any {\it specific} $n \geq 4$, 
but now the degree of the defining equation for $z$ is $2^{n-1}$, rather then $2^2=4$,
and the coefficients of that defining equation are polynomials in the $n(n+1)/2$ entries of the symmetric
$n \times n$ matrix $A$.

{\bf References}

[N1] Melvyn B. Nathanson, {\it Alternate minimization and doubly stochastic matrices}, \hfill \break
{\tt https://arxiv.org/abs/1812.11930} \quad .

[N2] Melvyn B. Nathanson, {\it Matrix scaling, explicit Sinkhorn limits, and Arithmetic}, \hfill \break
{\tt https://arxiv.org/abs/1902.04544} \quad .

[W] Avi Wigderson, {\it Alternate minimization and scaling algorithms: theory, applications, connections},
(The first of a 3-lecture series (the AMS Colloquium lectures)
given at the Joint Mathematics Meetings, San Diego, California - January 10, 2018). \hfill \break
Abstract, slides, and video available from {\tt https://www.math.ias.edu/avi/talks} \quad .

\bigskip
\hrule
\bigskip
Shalosh B. Ekhad and Doron Zeilberger, Department of Mathematics, Rutgers University (New Brunswick), Hill Center-Busch Campus, 110 Frelinghuysen
Rd., Piscataway, NJ 08854-8019, USA. \hfill\break
Email: {\tt [DoronZeil, ShaloshBEkhad] at gmail dot com}   \quad .
\bigskip
\hrule
\bigskip
Exclusively published in the Personal Journal of Shalosh B.  Ekhad and Doron Zeilberger and arxiv.org \quad .
\bigskip
 Written: Feb. 27, 2019.
\end